\documentclass[10pt,leqno]{amsart}
\usepackage{graphicx}
\usepackage{indentfirst,csquotes}

\usepackage{amssymb,amsthm,amsmath}
\usepackage{xcolor,paralist,hyperref,titlesec,fancyhdr,etoolbox}
\newtheorem{theorem}{Theorem}[]

\newtheorem{lemma}[theorem]{Lemma}
\newtheorem{proposition}[theorem]{Proposition}

\titleformat{\section}[display]{\normalfont\huge\bfseries\centering}{\centering\chaptertitlename\thechapter}{10pt}{\Large}
\titlespacing*{\section}{0pt}{0ex}{0ex}

\hypersetup{ colorlinks=true, linkcolor=black, filecolor=black, urlcolor=black }

\begin{document}
\title{{On synchronized coded systems}} 
\author[Wolfgang Krieger]{Wolfgang Krieger}
\date{\today}
\maketitle

\let\thefootnote\relax

\begin{abstract}
We introduce a class of codes with overlapping code words, that we call SPO-codes. The SPO-codes are related to the Markov codes that were introduced in: G. Keller,  J. Combinatorial Theory 56, (1991), pp.\ 75--83. The process of generating a coded system from a code extends to SPO-codes.  We describe a family of intrinsically ergodic synchronized SPO-coded systems that is closed under topological conjugacy. We construct synchronized subshifts with salient structural features by means of SPO-codes. We construct SPO-coded systems that are not semisynchronized.
\end{abstract} 

\bigskip


$\,$

$\,$
\begin{center}
1. Introduction
\end{center}

Let $\Sigma$ be a finite alphabet, and let $S_\Sigma$ denote the left shift on 
$\Sigma^\Bbb Z$. Given an $S_\Sigma$-invariant subset $X$ of $\Sigma^\Bbb Z$ the restriction $S_\Sigma$ of $S_\Sigma$ to $X$ is called a subshift.
 For basic notions of the theory of subshifts we refer to \cite {Ki} or \cite {LM}. A word is said to be admissible for a subshift $X\subset \Sigma^\Bbb Z$ if it appears on a point of $X$. We denote the set of admissible words of the subshift $X$ by $\mathcal L(X)$.
 
Coded systems were introduced in  \cite{BH}. Given the finite alphabet $\Sigma$, a set of words in the symbols of $\Sigma$ is called a code. The set of points in $\Sigma^\Bbb Z$ that carry a bi-infinite concatenation of words in $\mathcal C$ is called the concatenation set of the code.  A subshift that is the closure of the concatenation set of a code is called a coded system.
 
Before indicating the content of the paper we introduce terminology and notation.
Given a subshift $X \subset \Sigma^\Bbb Z$, and $i^-, i^+ \in \Bbb Z, i^-\leq i^+$ we
use the notation
$$
x_{[i^-, i^+]} = (x_i)_{i^- \leq i \leq i^+},\quad \quad \quad \quad
( x \in X),
$$
and the notation
$$
X_{[i^-, i^+]}   = \{x_{[i^-, i^+]}  : x \in X\}.
$$
We use similar notations in the case that indices range in left infinite or right infinite intervals. 
We identify the blocks in $X_{[1, l)}, l \in \Bbb N,$ with the words that they carry. 
We use the same notation for a finite or a left- or right-infinite block  
and for the finite or a left- or right-infinite word that the block carries. 
Concatenations of left-infinite, finite and right-infinite blocks and words we denote by juxtaposition.

Given a subshift we use notations like $x^{\langle - \rangle}$  $(x^{\langle + \rangle})$  for the points in $X_{(-\infty, i]} ( X_{[ i, \infty)} ),$ $i \in \Bbb Z$.
The right context of the finite block resp. finite word $a$ is denoted by 
\begin{align*}
&\Gamma^{\langle + \rangle}_\infty(a) = \{x^{\langle + \rangle} \in X_{(j, \infty)}: 
a x^{\langle + \rangle}\in X_{[i, \infty)}\}, \qquad 
(a \in  x_{[i,j]}, i,j \in \Bbb Z, i \leq j),
\\
&\Gamma^{\langle + \rangle}(a) = \{b \in \mathcal L(X): ab \in \mathcal L(X) \}, 
\qquad \qquad \qquad (a\in \mathcal L(X)).
\end{align*}
The notation $\Gamma^{\langle - \rangle}$ has the symmetric meaning.
We also set
\begin{align*}
\omega^{\langle + \rangle}(a) = \bigcap_{x^{\langle - \rangle}\in 
\Gamma_\infty^{\langle - \rangle}(a)} \{ x^{\langle + \rangle} \in \Gamma_\infty^{\langle + \rangle}&(a):
  x^{\langle - \rangle} a x^{\langle + \rangle} \in X\}, 
\\
 &(a \in  x_{[i,j]}, \  i,j \in \Bbb Z, i \leq j).
\end{align*}
The notation $\omega^{\langle - \rangle}$ has the symmetric meaning.

Given a subshift $X\subset \Sigma^\Bbb Z$ a point 
$ x^{\langle - \rangle} \in X_{(-\infty, i]}(x^{\langle + \rangle} \in X_{(-\infty, i]})$,$i \in \Bbb Z, $ is said to be transitive if for all $j < i$ ($j > i$) all words in $\mathcal L(X)$ appear in 
$x^{\langle - \rangle} _{(-\infty, j]} $($x^{\langle + \rangle}_{[j, \infty)} $).

Synchronized subshifts were introduced in \cite{BH}. A word $c \in \mathcal L(X)$ is synchronizing if one has for $u^{\langle - \rangle} \in \Gamma^{\langle - \rangle}(c)$ and 
$u^{\langle + \rangle} \in \Gamma^{\langle + \rangle}(c)$ that   
$u^{\langle - \rangle}cu^{\langle +\rangle}\in X$. A topologically transitive subshift is said to be synchronized if it has a synchronizing word. The set of synchronizing words of a synchronized subshift $X$ we denote by $\mathcal L_{synchro}(X)$. A word in 
$\mathcal L(X)$ that has a subword in $\mathcal L_{synchro}(X)$ is in $\mathcal L_{synchro}(X)$.

Semisynchronizing subshifts were introduced in \cite{Kr1}. A word $a$ is said to be semisynchronizing if $\omega^-(a)$ contains a transitive point. 
 A subshift is called semisynchronizing if it has a semisynchronizing word.

Hyposynchronized subshifts were introduced in \cite{Kr2}. A word $a \in \mathcal L(X)$ is hyposynchronizing iif $\omega^+(a)$ contains a transitive point,
and if there exists a point
$$
x^{\langle - \rangle} \in \Gamma^{\langle - \rangle}(a)
$$
such that
$$
\omega^{\langle + \rangle}(x^{\langle - \rangle}_{( - k, 0]} a)=
 \omega^{\langle + \rangle}(a),    \qquad (k \in \Bbb N).
$$ 
A subshift is hyposynchronized if it has a hyposynchronizing word. The inverse of a hyposynchronized subshift is semisynchronized.

For a subshift $X\subset \Sigma^\Bbb Z$ we denote by
$\mathcal L_{Markov}(X)$  the set of words $a \in \mathcal L(X) $ such that
$$
\# \{  \Gamma^+ (ba): b \in \Gamma^-(a) \}= \infty .
$$
The subshift that has $\mathcal L_{Markov}(X)$ as its set of admissibble words was introduced in \cite{T} and named the Markov boundary of $X$. We denote the Markov boundary of $X$ by $\partial_{Markov}(X)$.

The content of this paper is as follows.
Section 2 contains preliminaries. In particular it contains a lemma on synchronizing words and a lemma on semisynchronizing words. 

In Section 3 we  introduce a concatenation with overlap that we denote by $\circledast$ and suffix-prefix-codes with overlap (SPO-codes). The words in these codes come with a distinguished suffix and a distinguished prefix. In the concatenations that arise from such a SPO-code the distinguished suffix of a word and the distinguished prefix of the following word are the same word and this word appears only once in the concatenation. This yields a method of generating a subshift from a SPO-code that extends the method of generating a coded system from a code. Following the approach in \cite[Theorem 3]{KSW}, we give a sufficient condition on a SPO-code to generate an intrinsically ergodic subshift. We explain the relationship of the SPO-codes  to the Markov codes that were introduced in \cite{Ke}. 

In Section 4 we show that to a synchronized subshift $X\subset  \Sigma^\Bbb Z$ there is associated a PSO-code  $\mathcal C(X)$. This SPO-code  $\mathcal C(X)$ is canonical in the sense that its concatenation set is invariantly attached to $\mathcal C(X)$. As an application of the condition of Section 3 we give a  sufficient condition on a synchronizing subshift to be intrinsically ergodic. In the situation that we consider the set of synchronizing subshifts that satisfy this condition is closed under topological conjugacy. 

In Section 5 we construct by means of SPO-codes synchronized shifts with salient structural properties. We construct SPO-codes $\mathcal C$ such that 
$\partial_{Markov}(X(\mathcal C))$ consists of a finite number of orbits, and SPO-codes $\mathcal C$ such that 
$\partial_{Markov}(X(\mathcal C))$ consists of a countable infinity of orbits.

In Section 6 we describe  coded systems that are non-synchronized in the sense that they are neither semisynchronized nor hyposynchronized. These coded systems and their variations can be considered as structurally distant from the Dyck shifts.


\begin{center}
2. Preliminaries
\end{center}

We introduce more notation and terminology. Consider a finite alphabet $\Sigma$.
For $\sigma \in \Sigma$ we denote the word $a$ of length $n \in \Bbb N$ that is given by 
$a_i = \sigma, 1 \leq i \leq n,$ by $\sigma^n$.
The length of a finite word $a$ we denote by $\ell(a)$.
Beyond subshifts we consider more generally restriction $S_X$  of the shift $S_\Sigma$ to shift invariant Borel subsets of $\Sigma^\Bbb Z$. We call such restrictions Borel shifts. In the case that a Borel shift possesses  invariant probability measures its topological entropy $h(X)$ is defined as the supremum of the measure theoretic entropies of its invariant probability measures.

We recall that, given finite alphabets $\Sigma$ and $\bar \Sigma$, subshifts 
$X\subset  \Sigma^\Bbb Z  , \bar X \subset \bar \Sigma^\Bbb Z $, and a shift commuting continuous map
$$
\varphi: X \to \bar X,
$$
one has some $L \in \Bbb N$ such that $[-L. L]$ is the coding window of a block map
 $$
 \Phi: \Sigma^{[-L. L]} \to \bar \Sigma
  $$
that implements $\varphi$, in the sense that
$$
\varphi ( (x_i)_{i\in \Bbb Z}  ) = \Phi (x_{[i-L. L+i]})_{i\in \Bbb Z},  \quad \qquad \qquad
((x_i)_{i\in \Bbb Z} \in X).
$$
Given a coding window $[-L, L]$ and a bock map 
$$
 \Phi: \Sigma^{[-L. L]} \to \bar \Sigma
$$
we use the notation
$$
\Phi(x^{\langle - \rangle}) = (x^{\langle - \rangle}_{[i-L, i + L]})_{-\infty < i < j-L}, \qquad 
(x^{\langle - \rangle} \in X_{(-\infty, j}, j \in \Bbb Z).
$$

That the family of synchronized shifts is closed under topological conjugacy follows from the following lemma that
we include for completeness.

\begin{lemma}
Let $X \subset \Sigma^\Bbb Z,\bar X \subset \bar\Sigma^\Bbb Z$ be topologically conjugate subshifts and let
$$
\varphi: X \to \bar X
$$
be a topological conjugacy . Let $[-L, L]$ be a coding window for $\varphi$ and for
$\varphi^-$. 
Let $x \in X$ such that
\begin{align}
x_{[-L, L]} \in \mathcal L_{synchro},
\end{align}
and set $ \bar x = \varphi (x)$. Then 
\begin{align}
\bar x_{[-2L, 2L]}  \in \mathcal L_{synchro}(\bar X).
\end{align}
\end{lemma}
\begin{proof}
Let $\varphi$ be implemented by a block map
$$
 \Phi: \Sigma^{[-L. L]} \to \bar \Sigma
$$
and let $\varphi^-$ be implemented by a block map
$$
\bar \Phi: \bar\Sigma^{[-L. L]} \to  \Sigma.
$$
Let 
$$
 \bar y^{\langle - \rangle} \in \Gamma^{\langle - \rangle}(\bar x_{[-2L, 2L]}),
 $$
and let $y^{\langle - \rangle} \in X_{(-\infty, -L)}$ be given by
$$
y^{\langle - \rangle} = \bar \Phi (\bar y^{\langle - \rangle}, \bar x_{[-2L, 0)}  ).
$$
One has that 
$$
\bar\Phi (\bar y^{\langle - \rangle}, \bar x_{[-2L, 2L)}) = (y^{\langle - \rangle} ,  x_{[-L, L)}   ).
$$
Therefore
$$
y^{\langle - \rangle} \in \Gamma ^{\langle - \rangle}( x_{[-L, L)}).
$$

Symmetrically, let
$$
 \bar y^{\langle + \rangle} \in \Gamma^{\langle + \rangle}(\bar x_{[-2L, 2L]}),
 $$
and set
$$
y^{\langle + \rangle} = \bar \Phi (\bar y^{\langle + \rangle}, \bar x_{, 2L 0)}),
$$
and have that
$$
y^{\langle + \rangle} \in \Gamma ^{\langle + \rangle}( x_{[-L, L)}).
$$
It follows by (1) that
$$
(y^{\langle - \rangle} ,   x_{[-L, L)},  y^{\langle + \rangle}) \in X.
$$
One checks that
$$
\varphi (y^{\langle - \rangle} ,   x_{[-L, L)},  y^{\langle + \rangle})) =
(\bar y^{\langle - \rangle} ,   x_{[-2L, 2L)},  \bar y^{\langle + \rangle}
\in \bar X.
$$
This proves (2).
\end{proof}

For semisynchronizing words there is the following lemma.

\begin{lemma}
Let $X \subset \Sigma^\Bbb Z$ be a semisynchronizing subshift and let $a\in \mathcal L(X)$ be a semisynchronizing word of $X$. Let 
\begin{align}
c \in \Gamma^{\langle - \rangle}(a).
\end{align}
Then the word $ac$ is semisynchronizing.
\end{lemma}
\begin{proof}
Let
\begin{align}
x^{\langle - \rangle} \in \omega^{\langle - \rangle}(a)
\end{align}
be transitive. One shows that
\begin{align}
x^{\langle - \rangle} \in \omega^{\langle - \rangle}(ac).
\end{align}
For this one has to show for
\begin{align}
y^{\langle + \rangle} \in \Gamma_\infty^{\langle + \rangle}(ac),
\end{align}
that
\begin{align}
(x^{\langle - \rangle}, ac,y^{\langle + \rangle}) \in X.
\end{align}
One has from (3)
 and (6)
that
$$
cy^{\langle + \rangle} \in \Gamma_\infty^{\langle + \rangle}(a).
$$ 
From this one sees that (7)
 follows from (4).
\end{proof}

\begin{center}
3. Markov codes and SPO-codes
\end{center}

\noindent
3.1 $\bold{Markov \ codes.}$
Let there be given a finite alphabet $\Sigma$. 
A Markov code 
$^1$\footnote{$^1$ We consider codes that are Markov codes in the sense of Keller \cite{Ke}}
 is given by a code $\mathcal D \subset \Sigma^\Bbb Z$ as state space and an irreducible 0 - 1 
 transition$^2$ \footnote{$^2$For the countable state Bernoulli case see \cite{KSW}.} matrix 
$T_\mathcal D(d, d^\prime),d, d^\prime \in \mathcal D$.
The concatenation set
 $C(\mathcal D)$
of a Markov code 
$\mathcal D$ is the set of points 
$
x \in \Sigma^\Bbb Z
$ 
such that there exists a sequence of indices
$
i_k, k \in \Bbb Z,
$
$$
i_{k-1}< i_k, \qquad (k \in \Bbb Z,)
$$
such that
\begin{align}
&x_{(i_{k-1}, i_k]} \in \mathcal D,     
\\
&T_\mathcal D(x_{(i_{k-1}, i_k]}, x_{(i_{k}, i_{k+1}]}   )    = 1, \qquad  \qquad (k \in \Bbb Z). \notag
\end{align}
The Markov code is said to be unambiguous if for all points x in its concatenation set the sequence of indices$
i_k, k \in \Bbb Z,
$
that is assumed normalized,
$$
i_0 \leq 0 < i_1,
$$
 such that (8) holds, is unique. The closure of the concatenation set is the subshift$X(\mathcal D)$ generated by the Markov code.

\noindent
3.2. $\bold{SPO-codes.}$
For a finite alphabet $\Sigma$  
let there be given a bifix code $\mathcal F \subset \Sigma^+$. Denote by 
$\mathcal C_\mathcal F$ the set of words  $c \in \Sigma^+$ that have a proper prefix 
$f^-(c) \in \mathcal F$ and also a proper suffix $f^+(c) \in \mathcal F$. For a word 
$c \in \mathcal F $ we denote by $\mathring{c}$ the word that is obtained by removing from $c$ the suffix $f^+(c)$.

For
$a, b \in \mathcal C_\mathcal F$ we define the concatenation with overlap $a \circledast b$ by setting
$$
a \circledast b = \mathring{a} b.
$$
The word
$a \circledast b $ is in $\mathcal C_\mathcal F$:
$$
f^-(a \circledast b) = f^- (a), \quad f^+(a \circledast b) = f^+ (b).
$$
For $Q > 2$ and for
$$
c_q \in \mathcal C_\mathcal F, \quad \quad  1 \leq q \leq Q,
$$
we set
$$
\prod^\circledast_{1\leq q \leq Q}c_q =( \prod^\circledast_{1\leq q < Q}c_q ) \circledast c_Q.
$$

We define an SPO-code as a code $\mathcal C$ that is contained in 
$\mathcal C_\mathcal F$ for some bifix code $\mathcal F$. The concatenation set $C(\mathcal C)$ of an
SPO-code  $\mathcal C$ is the set of $x\in \Sigma^\Bbb Z$ such that there exists a sequence of indices
$
(j_j,i_k) \in \Bbb Z^2
$
such that 
\begin{align}
i_{k-1}  < j_k \leq i_k, \qquad (k \in \Bbb Z),
\end{align}
and such that 
\begin{align}
x_{(i_{k-1}.i_k ]} \in \mathcal C, \ \  x_{(j_{k}.i_k ]} =  f^+(x_{(i_{k-1}.i_k ]}  )  = 
f^-(x_{(i_k, i_{k+1}  ] }),    \qquad  (k \in \Bbb Z).
\end{align}
An SPO-code $\mathcal C$ is said to be unambiguous if for all points $x$ in its concatenation set the sequence. of indices  $(j_j,i_k) \in \Bbb Z^2,$
that is assumed normalized,
$$
i_0\leq0 < i_1,
$$
such that (8) and (9)  hold, is unique. The closure of $C(\mathcal C)$, to be denoted by $X(\mathcal C)$, is the subshift that is generated by $\mathcal C$.   

Given a bifix code $\mathcal F$ we set
$$
\mathcal C_\mathcal F^\bullet = \{c \in \mathcal C_\mathcal F^\circ: 
\ell(c) \leq \ell( f^-(c)) +\ell( f^+(c))\}.
$$
\begin{proposition}
Given a bifix code $\mathcal F$,  let  $\mathcal C  \subset \mathcal C_\mathcal F^\circ$ be an SPO-code
 with irreducible transition matrix, such that
$$
\mathcal C \cap  \mathcal C_\mathcal F^\bullet \neq \emptyset,
$$
and
\begin{align}
C(\mathcal C)) > h(X(\mathcal C) \setminus C(\mathcal C)),
\end{align}
and such that $X(\mathcal C)$ has a measure of maximal entropy of full support. Then $X(\mathcal C)$ 
is intrinsically ergodic.
\end{proposition}
\begin{proof}
Given a point $x \in C(\mathcal C)$ denote by
$$
(j_k(x),   i_k(x) )_{k \in \Bbb Z}
$$
the unique sequence 
$$
(j_k(x),   i_k(x) )_{k \in \Bbb Z} \in \Bbb Z^2,
$$
such that
$$
i_{k-1} < j_k \leq i_k, \qquad (k \in \Bbb Z),
$$
with the normalization
$$
i_0 \leq 0 < i_1,
$$
and such that
$$
x_{(i_{k-1}, i_k]}  \in \mathcal C,\quad x_{({j_k}, i_k]} = f^+(x_{(i_{k-1}, i_k]} ), \qquad 
(k \in \Bbb Z).
$$
We set
$$
\mathcal K(x) = \{k \in \Bbb Z:   x_{(i_{k-1}, i_k]}  \in \mathcal C_\mathcal F^\bullet  \},
\qquad \qquad (x \in C(\mathcal C)).
$$

From the unambiguous SPO-code $\mathcal C$ with irreducible transition matrix 
$T_{C}$ we obtain an unambiguous SPO-code  
$\widehat{\mathcal C}$
 with irreducible transition matrix $T_{\widehat C}$ by setting
 \begin{align}
 \widehat{\mathcal C}=
 \{(\prod^\circledast _{1\leq  q <Q}c_q)  \circledast c: c \in 
 \mathcal C_\mathcal F^\bullet, c_q \in 
 \mathcal C \setminus \mathcal C_\mathcal F^\bullet\, 1 \leq q < Q, Q \in \Bbb N\}.
  \end{align}
 Given a point $x \in C(\widehat{\mathcal C})$ denote by
$$
(\widehat{j}_k(x),   \widehat{i}_k(x) )_{k \in \Bbb Z}
$$
the unique sequence 
$$
(\widehat{j}_k(x),   \widehat{i}_k(x) )_{k \in \Bbb Z} \in \Bbb Z^2,
$$
such that
$$
\widehat{i}_{k-1} < \widehat{j}_k \leq \widehat{j}_k, \qquad (k \in \Bbb Z),
$$
that is assumed normalized,
$$
\widehat{i}_0 \leq 0 < \widehat{i}_1,
$$
and such that
$$
x_{(\widehat{i}_{k-1}, \widehat{i}_k]}  \in \widehat{\mathcal C},\quad 
x_{({j
\widehat{j}_k}, \widehat{i}_k]} = f^+(x_{(\widehat{i}_{k-1}, \widehat{i}_k]} ), 
\qquad (k \in \Bbb Z).
$$
We set
$$
\mathcal D_\mathcal C = \{\mathring{c}: c \in\widehat{\mathcal C} \}.
$$
By construction (see (12)) the words $d \in \mathcal D_\mathcal C $ have a prefix 
$f^-(d) \in \mathcal F$. We have obtained an unambiguous Markov code with state space
$\mathcal D_\mathcal C$ and transition matrix $T_{\mathcal D_\mathcal C}$ that is given by
$$
T_{\mathcal D_\mathcal C}(d, d^\prime) = 
\begin{cases}
1,\quad \text {if $df^-(d^\prime) \in  \widehat{\mathcal C}$},
\\
0, \quad \text {if $df^-(d^\prime)\not  \in  \widehat{\mathcal C} $}. 
\end{cases}
$$
 The concatenation set of the Markov code $\mathcal D_\mathcal C$ coincides with the concatenation set  of 
 $\widehat{\mathcal C} $.
 From the Markov code $\mathcal D_\mathcal C$ we obtain a topological Markov shift $Y(\mathcal C)$ with state space $\mathcal S$ given by
$$
\mathcal S =   \{(d, l): d \in \mathcal D_\mathcal C, 1 \leq l \leq \ell(d) \},
 $$
 and transition matrix $T_\mathcal S$ given by
 $$
 T_\mathcal S( (d, l), (. d, l+1)) = 1,  \quad d \in \mathcal D_\mathcal C, 1 \leq l < \ell (d),
 $$
 and
 $$
 T_\mathcal S( (d, \ell(d)), (d^\prime, 1) =
 \begin{cases}
  1, \quad  \text {if $T_{\mathcal D_\mathcal C}(d , d^\prime) = 1$},
 \\
  0, \quad  \text {if $T_{\mathcal D_\mathcal C}(d , d^\prime) = 0$}.
 \end{cases}
 $$
 As a dynamical system $C(\mathcal D_\mathcal C)$ is Borel conjugate to $Y(\mathcal C)$. A Borel conjugacy 
 $$
 \varphi: C(\mathcal D_\mathcal C) \to Y(\mathcal C)
 $$
 maps a point $x \in C(\mathcal D_\mathcal C) $ to the point $(s_i)_{i\in \Bbb Z}$ that 
 is determined by
 $$
 s_0(x) = (x_{(\widehat{i}_{0}, \widehat{i}_1]}   , -\widehat{i}_0).
 $$
 
 One has that
 $$
C( \widehat{\mathcal C}) \subset C( \mathcal C).
 $$
 The set of points $x \in  C( \mathcal C)$ such that $\mathcal K(x)$ is not bounded from above nor bounded from below is equal to the set $C( \widehat{\mathcal C})$. It follows by
 (11) that a measure of maximal entropy with full support of $X(\mathcal C)$ assigns measure one to $C( \widehat{\mathcal C})$, which is Borel conjugate to the countable state Markov shift with state space $\mathcal S$ and transititon matric $T_{\mathcal S}$. The hypothesis of the proposition therefore implies that this topological Markov shift has a measure of maximal entropy and therefore is positively recurrent, which means that it is intrinsically ergodic (\cite[Lemma  7.2.18]{Ki}). The intrinsic ergodicity of $X(\mathcal C)$  follows.
\end{proof}

In imposing the hypothesis (11) we have followed the approach of \cite[Theorem 3]{KSW}.

\begin{center}
4. Synchronization
\end{center}

\noindent
4.1. $ \bold{The  \ SP0-code \ of \
a \ synchronized \  subshift}.$
Given a synchronized subshift $X \subset \Sigma^\Bbb Z$ we set
\begin{align*}
\mathcal J_i(x) =& \{j \in (-\infty, i]: x_{[j, i]} \in \mathcal L_{synchro}(X)\}, \qquad 
(x\in X, i \in \Bbb Z),
\\
 &B^\circ(X) = \bigcap_{i\in \Bbb Z}\{x \in X: \mathcal J_i(x) \neq \emptyset \},
\\
J_i(x) = &\max \thinspace (\mathcal J_i(x)),\qquad \qquad \ \  \ \qquad
  (x\in X, i \in \Bbb Z, x \in B^\circ(X)).
\end{align*}
\begin{lemma}
For a synchronized subshift $X \subset \Sigma^\Bbb Z$ one has that
$$
J_{i+1}(x) \geq J_i(x), \qquad  \qquad (x\in X, i \in \Bbb Z, x \in B(X)).
$$
\end{lemma}
 \begin{proof}
 The lemma follows from
 $$
 x_{[J_i, i+1]}\in \mathcal L_{synchro}(X), \qquad \qquad (x\in X, i \in \Bbb Z). \qed
 $$
\renewcommand{\qedsymbol}{}
\end{proof}
 We denote
 $$
B(X) =  \{x\in B^\circ(X): \lim_{i \to \infty} J_i(x) = \infty \}.
$$

One notes that the set $B(X)$ is Borel.

\begin{lemma}
For an synchronized subshift $X$ that has an measure of full support $B(X)$ is not empty.
\end{lemma}
\begin{proof}
A measure of full support of the synchronized subshift $X$ assigns positive measure to $B(X)$.
 \end{proof}
 The next lemma shows that the set $B(X)$ is invariantly attached to $X$.
\begin{lemma}
Let $X \subset \Sigma^\Bbb Z,\bar X \subset \bar\Sigma^\Bbb Z$ be topologically conjugate subshifts and let
$$
\varphi: \bar X \to X
$$
be a topological conjugacy . Let $[-L, L]$ be a coding window for $\varphi$ and for
$\varphi^-$.Let $x \in X$ such that
$$
\lim_{i \to \infty} J_i(x) = \infty,
$$
and set $ \bar x = \varphi ^{-1}$. Then
$$
\lim_{i \to \infty}J_i(\bar x)= \infty.
$$
 \end{lemma}
 \begin{proof}
 Let $\bar k \in \Bbb Z$, and let $i \in \Bbb Z$ such that
 $$
 J_i(x) > \bar k + L.
 $$
 then
 $$
 x_{(J_i(x), i]}\in \mathcal L_{synchro}(X).
 $$
 By Lemma 1 then
 $$
 \bar x_{(J_i(\bar x) - L, i +L]}\in \mathcal L_{synchro}(\bar X).
 $$
 This means that
 $$
 J_{i+L}(\bar x) \geq \bar k.   
 $$ 
 \end{proof}
  
For $x \in  B(X)$ we denote by $\mathcal I(x)$ the set of $i \in \Bbb Z$ such that 
$$
J_i (x)> J_{i-1}(x).
$$
and we  enumerate $\mathcal I(x)$,
\begin{align}
\mathcal I(x) = \{I_i(x): i \in \Bbb Z \},
\end{align}
with the normalization
\begin{align}
I_0(x) \leq 0 < I_1(x), \qquad  (x \in B(X)).
\end{align}

We associate to the synchronizing subshift $X\subset \Sigma^\Bbb Z $ an SPO-code 
$\mathcal C(X)$.
The prefixes (suffixes) of the code $\mathcal C(X)$ arise as the words in
$$
\{x_{(J_i(x) , i]}:i \in \mathcal I(x), x \in B(X)\},
$$
and the words in $\mathcal C(X)$ arise as the words in
$$
\{x_{(I_i(x) , I_{i+1}(x)]}: i \in \Bbb Z, x \in B(X)\}.
$$
The  set $\mathcal F_{\mathcal C(X)} $ of prefixes (suffixes) of the code $\mathcal C(x)$ contains precisely the synchronizing words 
of $X$ that have no proper subword that is synchronizing. The words in the code 
$\mathcal C(X)$ are precisely the synchronizing words of $X$ that have a proper prefix and a proper suffix in $\mathcal F_{\mathcal C(X)} $ and that have no other subword that is synchronizing.

The set $B(X)$ is equal to $C(\mathcal C(X))$. As the indices  in  (4.4), (4.5) are uniquely determined by the point $x$, the SPO-code $\mathcal C(X)$ is unambiguous.

We formulate a Condition (H) on a synchronizing subshift $X$: 
\begin{align*}
\sup_{c \in \mathcal C(X)}\{\ell(c)  - \ell(f^- (c)) - \ell( f^+ (c))\} = \infty.  \tag H
\end{align*}

\begin{lemma}
The condition H is an invariant of topological conjugacy of synchronizing shifts.
\end {lemma}
\begin{proof}
Apply Lemma 1.
\end{proof}

\begin{theorem}
Let $X$ be a 
topologically transitive synchronized subshift that has a measure of maximal entropy of full support, that satisfies Condition H, and is such that
$$
h(B(X))  > h(X \setminus B(X)).
$$
Then $X$ is intrinsically ergodic.
\end{theorem}
\begin{proof}
Apply Proposition 3.
\end{proof}

 \begin{center}
 5. A class of synchronized shifts
 \end{center}

We give an example  such that 
$\partial_{Markov}(X(\mathcal C))$ contains finitely many orbits and an example of an SPO-code $\mathcal C$ such that 
$\partial_{Markov}(X(\mathcal C))$ contains a countable infinity of orbits. The SPO-coded shifts of the examples satisfy Condition H, and are therefore intrinsically ergodic by Theorem 8.

For the construction of the SPO-codes we use the alphabet 
$\{\delta, \gamma\}$ and a finite alphabet $\Sigma$ that is disjoint from 
$\{\delta, \gamma\}$. 
The words $\gamma\delta^n\gamma, n \in \Bbb N,$ will serve as the distinguished suffixes and prefixes of the SPO-codes.

\noindent
5.1. $\bold {Example  \ 1.}$ 
We construct an SPO-code $\mathcal C$ such that
$\partial_{Markov}(X(\mathcal C))$ is a finite set of orbits. For the construction we choose a finite set 
$P$ of periodic points of $(\Sigma^\Bbb Z, S_{\Sigma^\Bbb Z})$.

Set
$$
R = \prod_{p \in  P} \text{period of} (p).
$$
 The code is given by
 
\begin{align*}
\mathcal C = &\{\gamma \delta \gamma p_{[0   , 2qR   )} \gamma \delta \gamma:
p \in P, q > 1  \} \cup
\\
&\{\gamma \delta^n \gamma p_{[0   , 2qR   )} \gamma \delta^{n+1} \gamma : p \in P,   
n \in \Bbb N, q \geq n\} \cup
\\
&\{\gamma \delta^n \gamma p_{[0   , 2qR)} \delta^{n-1} :  p \in P, n > 1,  q \geq n \} .
\end{align*}

One has that
$$
\Gamma (\delta^{n+1} \gamma p_{[0, R)})
\begin{cases}
\ni p_{[0,n R)}\gamma,
\\
  \not \ni p_{[0,(n +1)R)}\gamma, 
 \end{cases}
  \qquad \qquad (p \in  P, n \in \Bbb N).
$$

It follows that $\partial_{Markov} (X(\mathcal C))$ contains the points in $P$, the fixed point $(x_i)_{i\in \Bbb Z}$ that is given by $x_i = \delta, i \in \Bbb Z,$ and the orbits that contain the points 
$x^{p, -},x^{(p, +)}, p \in  P,$
that are given by, 
$$
x^{p, -}_{(-\infty, 0)} = p _{(-\infty, 0)},   \qquad \qquad x^{p, +}_{[0, \infty)} = p_{[0, \infty)},
$$
$$
\qquad x^{p, -}_{i} =
\begin{cases}
\gamma, \text{if $i = 1,  $}
\\
\delta, \text{if $i > 1 , $}
\end{cases}
\quad \qquad \qquad
x^{p, +}_{i} =
\begin{cases}
\gamma, \text{if $i = 1,  $}
\\
\delta, \text{if $i > 1 , $}
\end{cases}
$$
and no other points.

\noindent
5.2. $\bold{Example  \ 2.}$
We construct an SPO-code $\mathcal C$ such that
$\partial_{Markov}(X(\mathcal C))$ is countable infinity of orbits. We denote by $p^{(k)}$ the point of period $k$ in $(\{0, 1\}^\Bbb Z,S_{ \{0, 1\}^\Bbb Z})$ that is given by
$$
p^{(k)}_i=
\begin{cases}
0, \text{if $0 \leq i < k -1 $,}
\\
1,\text{if $ i = k -1 $.}
\end{cases}
$$

Set
$$
R_k = 1 + \tfrac{1}{2}{k(k-1)}, \qquad (k \in \Bbb N).
$$

The code is given by
\begin{align*}
\mathcal C = &\{\gamma  \delta^k \gamma p^{(1)}_{[0   , 2k    )} \gamma\delta^{(R_k)}\gamma:
k \in \Bbb N \} \cup
\\
&\{\gamma  \delta^k \gamma p^{(m)}_{[0   , 2mk    ) } \gamma  \delta^{k + m -1} \gamma:
  1 < m \leq k, k \in \Bbb N \}.
\end{align*}

Consider the words 
$$
b^-(q,k,m) = \delta^q \gamma  p^{(k)}_{[0, 2mk)}, \qquad (q, k, m \in \Bbb N),
$$
that are in $\mathcal L (\mathcal C(X))$.
Set
$$
K^-_\circ =
\begin{cases}
k - q, \text{if $q<k,$}
\\
0, \text{if $q=k$.}
\end{cases}
$$

One has that
$$
\Gamma^+ (\gamma\delta^{K} b^-(m,k,q))
\begin{cases}
\ni p^{(m)}_{[0,2mk)}\gamma,
\\
  \not \ni p_{[0,2m(K+1)}\gamma, 
 \end{cases}
  \qquad \qquad (K \geq K_\circ).
$$

This implies that $b^-(q,k,m)\in \mathcal L_{Markov}(\mathcal C(X))$.

Consider the words 
$$
b^+(q,k,m) = \delta^q \gamma  p^{(k)}_{[0, 2km)}, \qquad (q, k, m \in \Bbb N),
$$
that are in $\mathcal L (\mathcal C(X))$.
Set
$$
M^-_\circ =
\begin{cases}
q - m, \text {if $q>m,$}
\\
0, \text {if $q\leq m$.}
\end{cases}
$$

One has that
$$
\Gamma^+(\gamma p^{(k)}_{[0,2kM)}b^+(q,k,m))
\begin{cases}
\ni\delta^{ R_{\max\{q, m\} + M - M_\circ} +\max\{q, m\} + M - M_\circ -1 - q}\gamma,
\\
\not\ni\delta^{ R_{\max\{q, m\} + M - M_\circ} +\max\{q, m\} + M - M_\circ  - q}, 
 \end{cases}(M \geq M_\circ).
$$

This implies that $b^+(q,k,m)\in \mathcal L_{Markov}(\mathcal C(X))$.

It follows that $\partial_{Markov}(X(\mathcal C))$ contains the points $p^{(k)}, k \in \Bbb N$, the  the fixed point $(x_i)_{i\in \Bbb Z}$ that is given by $x_i = \delta, i \in \Bbb Z,$  the fixed point $(y_i)_{i\in \Bbb Z}$ that is given by $y_i = 0, i \in \Bbb Z,$  and the orbits that contain the points $x^{k, -},x^{(k, +)}, k \in  \Bbb N,$ that are given by, 
$$
x^{(k, -)}_{(-\infty, 0)} = p^{(k)}_{(-\infty, 0)},   \qquad \qquad x^{p, +}_{[0, \infty)} = p^{(k)}_{[0, \infty)},
$$
$$
\qquad x^{p, -}_{i} =
\begin{cases}
\gamma, \text{if $i = 1,  $}
\\
\delta, \text{if $i > 1 , $}
\end{cases}
\quad \qquad \qquad
x^{p, +}_{i} =
\begin{cases}
\gamma, \text{if $i = 1,  $}
\\
\delta, \text{if $i > 1 , $}
\end{cases}
$$
and no other points.

\begin{center} 
6. A class of  coded systems
\end{center}

For the construction of  examples of subshifts that are not semisynchronizing we use the alphabet
$$
\Sigma = \{-1, 0, 1\}.
$$
For $k \in \Bbb Z_+, \alpha \in \{-1, 1\}$, we denote by $g^{\langle -, k, \alpha \rangle}$
the word of length $k+1$ that is given by
\begin{align*}
&g^{\langle -, k, \alpha \rangle}_{i} =  
\begin{cases}
0, \ \text{if $  1 \leq i \leq k,$}
\\
 \alpha, \ \text{if $ i = k+1,$}
\end{cases}
\end{align*}
and we denote by $g^{\langle +, k, \alpha \rangle}$
the word of length $k+1$ that is given by
\begin{align*}
&g^{\langle -, k, \alpha \rangle}_{i} =  
\begin{cases}
\alpha, \ \text{if $ i = 1,$}
\\
 0, \ \text{if $ 1 <i \leq k+1.$}
\end{cases}
\end{align*}
For $k \in \Bbb N, \alpha \in \{-1, 1\}$ we denote by $c^{\langle k, \alpha \rangle}$, the word of length $k + 2$ that is given by
$$
c^{\langle k, \alpha \rangle}_{i}=
\begin{cases}
\alpha, \ \text{if $ i =1.$}
\\
0, \ \text{if $ 1 <i \leq k+1.$}
\\
\alpha \ \text{if $ i = k+2.$}
\end{cases}
$$

Note that
\begin{align*}
\alpha &= g^{\langle -, 0, \alpha \rangle} = g^{\langle +, 0, \alpha \rangle}, \qquad \quad
(\alpha \in \{-1, 1\}).
\\
c^{\langle k, \alpha \rangle} &= g^{\langle +, k^+, \alpha \rangle} 
  g^{\langle -, k^-, \alpha \rangle},
\qquad  \ (\alpha \in \{-1, 1\}, k^+, k^- \in \Bbb Z_+, k^++ k^- > 0).
\end{align*}

We denote the coded system that is generated by the code
$$
\mathcal C =\{ c^{\langle k, \alpha \rangle}  :  k \in \Bbb N, \alpha \in \{-1, 1\}\}
$$
by $X_\mathcal C $. We set
\begin{align*}
&\mathcal D^{\langle 1 \rangle} = \{\alpha 0 \alpha: \alpha \in \{-1, 1\} \},
\\
&\mathcal D^{\langle  m + 1 \rangle} = \{\alpha (\bigcup_{1\leq l \leq m}
\mathcal D^{\langle  l \rangle}) ^+\alpha: \alpha \in \{-1, 1\} \}, \qquad (m \in \Bbb N).
\end{align*}

We obtain a subshift $X \subset \Sigma^\Bbb Z$ by excluding from $X_\mathcal C$ the words in
$$
\bigcup _{k\in \Bbb N, \alpha \in \{-1, 1\}}
\{c^{\langle k, \alpha^- \rangle}c^{\langle k+m, \alpha^+ \rangle},c^{\langle k, \alpha^- \rangle}c^{\langle k+m, \alpha^+ \rangle}: m > 1, \ \alpha^-, \alpha^+\in \{-1, 1\}  \}.
$$
and the words in
$
\bigcup_{m\in \Bbb N} \mathcal D^{\langle  m \rangle}.
$

Remark:
For a word $a \in \mathcal L(X)$ that does not have a word in $\mathcal C$ as a subword one can find a word $b \in \mathcal L(X)$ that has the word $a$ as a prefix and that has a word in $\mathcal C$ as a subword: We can distinguish four cases. 

For $\alpha \in \{-1, 1\}, k \in \Bbb Z_+,$ and
$$
a = g^{\langle +, k, \alpha \rangle}, 
$$ 
set $b$ equal to 
$$
g^{\langle +, k, \alpha \rangle} g^{\langle -, 0, \alpha \rangle} = 
c^{\langle k, \alpha^- \rangle}.
$$

For $(\alpha \in \{-1, 1\}, k \in \Bbb N),$ and
$$
a = g^{\langle -, k, \alpha \rangle}, 
$$
set $b$ equal to 
$$
g^{\langle -, k, \alpha \rangle}  c^{\langle k+1, \alpha \rangle},
$$ 
which is a subword of 
$$
g^{\langle +, 0, \alpha \rangle}  g^{\langle -, k, \alpha \rangle}  
c^{\langle k+1, \alpha \rangle} =c^{\langle k, \alpha \rangle}  
c^{\langle k+1, \alpha \rangle}.
$$

For $\alpha \in \{-1, 1\}, k^+ > k^-,$ and
$$
a = g^{\langle -, k^-, \alpha^- \rangle}   g^{\langle +, k^+, \alpha^+ \rangle},
 $$
set $b$ equal to 
$$
a g^{\langle -, 0, \alpha^+ \rangle} = g^{\langle -, k^., \alpha^- \rangle} 
c^{\langle k+, \alpha^+ \rangle}
$$ 
which is a subword of 
$$
g^{\langle +, k^+ .k^- -1, \alpha^- \rangle}  g^{\langle -, k^- \alpha^- \rangle}c^{\langle k^+, \alpha^+ \rangle}
=c^{\langle k^+ -1, \alpha^- \rangle} 
c^{\langle k^+, \alpha^+ \rangle}.
 $$
 
For $\alpha \in \{-1, 1\}, k^+ \leq k^-$, and
$$
a = g^{\langle -, k^-, \alpha^- \rangle}  g^{\langle -, k^+, \alpha^+ \rangle},
  $$
set $b$ equal to 
$$
a   g^{\langle -, k^- - k^+ +1, \alpha^+ \rangle} = g^{\langle -, k^-, \alpha^- \rangle} 
c^{\langle k^- + 1, \alpha^+ \rangle}
$$ 
which is a subword of 
$$
g^{\langle +, 0, \alpha^- \rangle}   g^{\langle -, k^-, \alpha^- \rangle} 
 c^{\langle k^- + 1, \alpha^+ \rangle} =
 c^{\langle k^-, \alpha^- \rangle} 
c^{\langle k^- + 1, \alpha^+ \rangle}.
 $$

\begin{lemma}
 Let $R \in \Bbb N$, and
$$
k_r \in \Bbb N,  \alpha_r \in \{-1, 1\},\qquad (1 \leq r \leq R).
$$
such that
$$
a =  \prod_{1 \leq r \leq R} 
c^{\langle k_r, \alpha_r \rangle} \in \mathcal L(X).
$$
Assume that
$$
Q  = \max\{k_r:1 \leq r \leq R\} - k_R > 0.
$$
Then there exist
$$
\alpha _q \in \{-1, 1\}, \qquad (1 \leq q \leq Q),
$$
such that
$$
a  \prod_{1 \leq q \leq Q} c^{\langle k_R + q, \alpha_q \rangle} \in \mathcal L(X).
$$
\end{lemma}
\begin{proof}
Set
\begin{align*}
&r_1 = \max \{r \in [1, R): k_r = k_R +1\} ,
\\
&r_{q+1} = \max \{r \in [1, r_q): k_r = k_q+1\}, \qquad (1 \leq q < Q), 
\end{align*}
and set
$$
\alpha_q = \alpha^-_{r_{q}}, \qquad (1 \leq q \leq Q). 
$$
\end{proof}
\begin{lemma}
Let $k \in \Bbb Z_+, \alpha \in \{-1, 1\}$. Also let $R \in \Bbb N$,
$$
k_r \in \Bbb N,  \alpha_r \in \{-1, 1\},\qquad (1 \leq r \leq R),
$$
such that 
$$
k_R \geq k_r, \qquad (1\leq r < R).
$$

(a) Assume that  
$$
a =  \prod_{1 \leq r \leq R} 
c^{\langle k_r, \alpha_r \rangle} \in \mathcal L(X).
$$
Then the word $a$ is not semisynchronizing.

(b) Assume that
$$
a = g^{\langle -, k, \alpha \rangle}  \prod_{1 \leq r \leq R} 
c^{\langle k_r, \alpha_r \rangle} \in \mathcal L(X)..
$$
Then the word $a$ is not semisynchronizing.
\end{lemma}
\begin{proof}
(a)
Let
$$
x^{\langle - \rangle} \in \omega^{\langle -  \rangle}(a)
$$
be transitive. Let $i_\circ$ be the maximal index in $(-\infty, 0]$ such that
$$
x^{\langle - \rangle} _{(i - k_R - 3 ,i]} \in \{c^{\langle k_R + 1, \alpha \rangle}:
 \alpha \in \{-1, 1 \}\}.
$$
Then
$$
y^{\langle + \rangle} = x^{\langle - \rangle} _{(i_\circ - k_R - 3 ,i_\circ]} 
\prod_{1 < q < \infty}  c^{\langle k_R + q, 1 \rangle}  \in \Gamma^{\langle + \rangle}(a),
$$
but
$$
( x^{\langle - \rangle}  ,  a , y^{\langle + \rangle} ) \not \in X.
$$

(b) Let
$$
x^{\langle - \rangle} \in \omega^{\langle -  \rangle}(a)
$$
be transitive. In the case that $ x^{\langle - \rangle}$ has the suffix 
$$
g^{\langle +, k_R - k + 1, \alpha \rangle}
$$
let $i_\circ$ be the maximal index in $(-\infty, 0]$ such that
$$
x^{\langle - \rangle} _{(i - k_R - 4 ,i]} \in \{c^{\langle k_R + 2, \alpha \rangle}:
 \alpha \in \{-1, 1 \}\}.
$$
Then
$$  
y^{\langle + \rangle} =   c^{\langle k_R + 1, \alpha^{-1} \rangle}  
x^{\langle - \rangle} _{(i - k_R - 4 ,i]}
\prod_{2 < q < \infty}  c^{\langle k_R + q, 1 \rangle}  \in \Gamma^{\langle + \rangle}(a)
    \in \Gamma^{\langle + \rangle}(a),
$$
but
$$
( x^{\langle - \rangle}  ,  a , y^{\langle + \rangle} ) \not \in X.
$$

Otherwise
let $i_\circ$ be the maximal index $i$ in $(-\infty, 0]$ such that
$$
x^{\langle - \rangle} _{(i - k_R - 3 ,i]} \in \{c^{\langle k_R + 1, \alpha \rangle}:
 \alpha \in \{-1, 1 \}\}.
$$
Then
$$
y^{\langle + \rangle} =   c^{\langle k_R + 1, x^{\langle - \rangle}_{i_\circ} \rangle} 
\prod_{1 < q < \infty}  c^{\langle k_R + q, 1 \rangle}  \in \Gamma^{\langle + \rangle}(a),
$$
but
$$
( x^{\langle - \rangle}  ,  a , y^{\langle + \rangle} ) \not \in X. 
$$
\end{proof} 
\begin{proposition}
The subshift $X$ is not semisynchronizing
\end{proposition}
\begin{proof}
The subshift $X$ cannot have a semisynchronizing word. This follows by Lemma 1 from the remark by an application of Lemma 9 and of Lemma 10. 
\end{proof}

The construction of the subshift $X$ is time-symmetric. It follows  that the subshift $X$ is not hyposynchronizing.

\par\noindent Wolfgang Krieger
\par\noindent Institute for  Mathematics, 
\par\noindent  University of Heidelberg,
\par\noindent Im Neuenheimer Feld 205, 
 \par\noindent 69120 Heidelberg,
 \par\noindent Germany
\par\noindent krieger@math.uni-heidelberg.de

\end{document}